\documentclass[reqno, A4paper]{amsart} 
\usepackage{amsmath, amssymb, amsfonts, amstext, amsthm, amscd}
\usepackage[mathscr]{euscript}
\usepackage{enumerate}
\usepackage{tikz-cd}
\usepackage[Symbolsmallscale]{upgreek}
\usepackage{mathrsfs}
\usepackage{eqlist}
\usepackage{array}
\usepackage{tikz,color,soul}
\usetikzlibrary{matrix,arrows}
\usepackage[english]{babel}
\usepackage{chngcntr}
\usetikzlibrary{arrows.meta, positioning,arrows}
\theoremstyle{plain}
\newtheorem{thm}{Theorem}[section]
\newtheorem{prop}[thm]{Proposition}
\newtheorem{lemma}[thm]{Lemma}
\newtheorem{cor}[thm]{Corollary}
\theoremstyle{defn}
\newtheorem{defn}[thm]{Definition}
\newtheorem{ex}[thm]{Example}
\newtheorem{rmk}[thm]{Remark}
\newtheorem{conj}[thm]{Conjecture}

\numberwithin{equation}{section}

\newcommand{\co}{\mathcal O}
\newcommand{\C}{\mathbb C}

\newcommand{\N}{\mathbb N}

\newcommand{\R}{\mathbb R}
\newcommand{\Z}{\mathbb Z}

\begin{document}

\title{Cohomology of Torus Manifold Bundles}

\author{Jyoti Dasgupta}
\address{Department of Mathematics, Indian Institute of Technology-Madras, Chennai, India}
\email{jdasgupta.maths@gmail.com}

\author{Bivas Khan}
\address{Department of Mathematics, Indian Institute of Technology-Madras, Chennai, India}
\email{bivaskhan10@gmail.com}

\author{V. Uma}
\address{Department of Mathematics, Indian Institute of Technology-Madras, Chennai, India}
\email{vuma@iitm.ac.in}


\begin{abstract} Let $X$ be a $2n$-dimensional torus manifold with a locally standard $T \cong \left( S^1 \right)^n$ action whose orbit space is a homology polytope. Smooth complete complex toric varieties and quasitoric manifolds are examples of torus manifolds. Consider a principal $T$-bundle $p : E \rightarrow B$ and let $\pi : E(X) \rightarrow B$ be the associated torus manifold bundle. We give a presentation of the singular cohomology ring of $E(X)$ as a $H^*(B)$-algebra and the topological $K$-ring of $E(X)$ as a $K^*(B)$-algebra with generators and relations. These generalize the results in \cite{masudapan} and \cite{param} when the base $B=pt$. These also extend the results in \cite{paramuma}, obtained in the case of a smooth projective toric variety, to any smooth complete toric variety.  \end{abstract}

\maketitle
\section{Introduction}

A torus manifold is a $2n$-dimensional manifold acted upon effectively by an $n$-dimensional compact torus with non-empty fixed point set. Smooth complete complex toric varieties and quasitoric manifolds are examples of torus manifolds. The notion of torus manifolds was introduced by A. Hattori and M. Masuda in \cite{hm}. In \cite{masudapan} M. Masuda and T. Panov studied relationships between the cohomological properties of torus manifolds and the combinatorics of their orbit spaces. The topological $K$-ring of the torus manifolds with locally standard action and orbit space a homology polotope was described by P. Sankaran in \cite{param}.

Let $p : E \rightarrow B$ be a principal bundle with fibre and structure group the complex algebraic torus $\mathbb{T} \cong (\C^*)^n$ over a topological space $B$. For a smooth projective $\mathbb{T}$-toric variety $X$, consider the toric bundle $\pi : E(X) \rightarrow B$, where $E(X)=E \times_{\mathbb{T}} X$, $\pi([e, x])=p(e)$. In \cite{paramuma}, the authors describe the singular cohomology ring of $E(X)$ as a $H^*(B)$-algebra. Furthermore, when $B$ is compact Hausdorff, they describe the topological $K$-ring of $E(X)$ as a $K^*(B)$-algebra. 

In this paper we consider $p: E \rightarrow B$ to be a principal bundle with fibre and structure group the compact torus $T \cong \left( S^1 \right)^n$.  We assume that $B$ has the homotopy type of a finite CW complex so that $H^*(B)$ and $K^*(B)$ are finitely generated abelian groups. Without loss of generality, we further assume that $B$ is compact and Hausdorff. Let $X$ be a $2n$-dimensional torus manifold with a locally standard action of $T$ such that the orbit space $Q:=X/T$ is a homology polytope. We call the associated bundle $\pi : E(X) \rightarrow B$ a torus manifold bundle, where $E(X)=E \times_T X$. In Theorem \ref{2} we give a presentation of the singular cohomology ring of $E(X)$ as a $H^*(B)$-algebra. A presentation of the topological $K$-ring $K^*(E(X))$ as a $K^*(B)$-algebra is obtained in Theorem \ref{6}. As an application, we describe the cohomology ring and $K$-ring of toric bundles for a smooth complete toric variety in Corollary \ref{15} extending the results in \cite{paramuma}.

The method of proof for Theorem \ref{2} exploits the known presentation of the cohomology ring \cite[Corollary $7.8$]{masudapan} when the base $B$ is a point. Applying the Leray-Hirsch theorem in cohomology we first prove that $H^*(E(X))$ is a free module over $H^*(B)$ of rank $\chi(X)$. Then we construct a surjective $H^*(B)$-algebra homomorphism from $R(B, \left( Q, \Lambda \right))$ (see Definition \ref{3}) to $H^*(E(X))$. Here $\Lambda$ denotes the characteristic map of the torus manifold (see Section 2).  To verify that this algebra homomorphism is injective, we recall from \cite{masudapan} that the equivariant cohomology ring which is isomorphic to the face ring of $Q$, is a free $H^*(BT)$-module of rank $\chi(X)$, where $BT$ denotes the classifying space of principal $T$-bundles. We then canonically extend the scalars of the face ring to $H^*(B)$ and use that it is a finitely generated abelian group to conclude injectivity.

Similarly, the method of proof for Theorem \ref{6} exploits the known presentation of the topological $K$-ring \cite[Theorem $5.3$]{param} when the base $B$ is a point. Applying the Leray-Hirsch theorem in $K$-theory we first prove that $K^*(E(X))$ is a free module over $K^*(B)$ of rank $\chi(X)$.  Then we construct a surjective $K^*(B)$-algebra homomorphism from $\mathcal{R}(B, \left( Q, \Lambda \right))$ (see Definition \ref{5}) to $K^*(E(X))$. Let $M:=\text{Hom}(T,{S}^1)$ denote the character lattice of $T$ and $RT:=\mathbb{Z}[\chi^u: u\in M]$ the ring of finite dimensional complex representations of $T$. In Proposition \ref{8} we show that the {\it $K$-theoretic face ring} of $Q$ denoted by $\mathcal{K}(Q)$ (see Definition \ref{14}) is a free $RT$-module of rank $\chi(X)$, using methods similar to \cite{vezzosi2003higher} and \cite{baggio2007equivariant} in the setting of smooth toric varieties. We then canonically extend the scalars of $\mathcal{K}(Q)$ to $K^*(B)$ and use that it is a finitely generated abelian group to conclude injectivity. In the case of a smooth complete toric variety, the $K$-theoretic face ring is in fact isomorphic to the algebraic $\mathbb{T}$-equivariant $K$-ring \cite[Theorem 6.4]{vezzosi2003higher}. The authors believe that the topological equivariant $K$-ring of any $T$-torus manifold is isomorphic to the $K$-theoretic face ring but could not find it in literature. We prove this statement for a quasitoric manifold in a parallel work \cite{dku}.

In Section 6 we consider a torus manifold $X$ with a locally standard action of $T$ such that $X/T$ is not necessarily a homology polytope but only a face-acyclic nice manifold with corners (see Section 2 for the definition). The equivariant cohomology ring as well as the ordinary cohomology ring of $X$ have been described by Masuda and Panov in \cite[Theorem 7.7, Corollary 7.8]{masudapan}. Let $E(X)\longrightarrow B$ be the bundle with fiber $X$ associated to the principal $T$-bundle over a topological space $B$ which is of the homotopy type of a finite CW complex. We generalize \cite[Corollary 7.8]{masudapan} to give a presentation $H^*(E(X))$ as a $H^*(B)$-algebra in Theorem \ref{cohombuntorus}. Similar to Theorem \ref{2} we prove this by using the Leray-Hirsch theorem and the known presentation of $H^*_{T}(X)$ as a $H^*(BT)$-algebra \cite[Theorem 7.7]{masudapan}. We finally conjecture a similar presentation for $K^*(E(X))$ as a $K^*(B)$-algebra. We note that difficulties arise in extending the result to this setting especially because the cohomology ring $H^*(X)$ is not generated in degree $2$.

\noindent {\bf Acknowledgements:} The authors are grateful to Prof. P. Sankaran for drawing our attention to this problem and for his valuable comments on the initial versions of this manuscript. The first and the second author thank the Council of Scientific and Industrial Research (CSIR) for their financial support.  The authors wish to thank the unknown referee for a careful reading of the manuscript and for very valuable comments and suggestions which led to improving the text. The final section has been added taking into account the referee's suggestions. The extension of Theorem 3.3 to Theorem 6.1 was also suggested by Prof. M. Masuda in a prior email correspondence. We are grateful to him for this.

\section{Notation and Preliminaries}
We recall some notation and preliminaries from \cite{masudapan} and \cite{param}. 

\subsection{Torus manifolds}

Let $T \cong \left( S^1 \right)^n$ denote the compact $n$-dimensional torus. A $2n$-dimensional closed connected orientable smooth manifold $X$ with an effective smooth action of $T$ such that the fixed point set $X^T$ is non-empty, is called a \emph{torus manifold}. Since $X$ is compact it follows that $X^T$ is finite (see \cite[Section 3.4]{gmukherjee},  \cite[Section 7.4]{bp}).  A codimension-two connected submanifold is called a \emph{characteristic submanifold} of $X$ if it is pointwise fixed by a circle subgroup of $T$. Since $X$ is compact, there are finitely many characteristic submanifolds, which we denote by $V_1, \ldots, V_d$. It can be shown that each $V_i$ is orientable. We say that $X$ is omnioriented, if an orientation is fixed for $X$ and for every characteristic submanifold $V_i$. We fix an omniorientation of $X$.
 
The $T$-action on the torus manifold $X$ is said to be locally standard if it has a covering by $T$-invariant open sets $U$ such that $U$ is weakly equivariantly diffeomorphic to an open subset $U'\subset \C^n$ invariant under the standard $T$-action on $\C^n$. The latter means that there is an automorphism $\theta :T \rightarrow T$ and a diffeomorphism ${g}: U \rightarrow U'$ such that ${g}(ty)=\theta(t) {g}(y)$ for all $t \in T$, $y \in U$. Let $Q:=X/T$ be the orbit space and let $\Upsilon : X \rightarrow Q$ be the projection map. If $X$ is locally standard, then $Q$ becomes a nice manifold with corners (see \cite[Section 5.1 p. 724]{masudapan} \cite[Definition 7.1.3]{bp}). We denote by $Q_i$ the image of $V_i$ under $\Upsilon$ for $i=1, \ldots, d$; these are the \emph{facets} or the codimension one faces of $Q$. A codimension-$k$ \emph{preface} is defined to be a non-empty intersection of $k$ facets for $k=1,\ldots, n$. The connected components of prefaces are called \emph{faces}. We regard $Q$ itself as a face of codimension zero.  We say that $Q$ is {\em face-acyclic} if all its faces are acyclic i.e. $\tilde{H}_{i}(F)=0,~\mbox{for every}~ i$, for each face $F$ of $Q$. We say that $Q$ is a \emph{homology polytope}, if $Q$ is face-acyclic and all its prefaces are faces. This is equivalent to saying that $Q$ is acyclic and all its prefaces are acyclic (in particular, connected). In this case the intersection of $r$ facets $Q_{i_1}, \ldots, Q_{i_r}$ is a codimension $r$ face $F$ of $Q$. Equivalently non-empty intersections of characteristic submanifolds are connected submanifolds of $X$. Unless otherwise specified, we shall assume henceforth that $X$ is a locally standard torus manifold with $Q$ a homology polytope. Note that, $H^{\ast}(X)$ is generated in degree two if and only if $X$ is locally standard and $Q$ is homology polytope (see \cite[Theorem $8.3$]{masudapan}).  For every characteristic submanifold $V_i$, there is a primitive element $v_i \in \text{Hom }(S^1, T) \cong \Z^n$ determined up to sign, whose image is the circle subgroup fixing $V_i$ pointwise. The sign of $v_i$ is determined by the omniorientation. Define the characteristic map $\Lambda : \{ Q_1, \ldots, Q_d\} \rightarrow \text{Hom}(S^1, T)$, such that $\Lambda(Q_i)=v_i$. The local standardness of $X$ implies that the characteristic map $\Lambda$ satisfies the following smooth condition: if $Q_{i_1}\cap \cdots \cap Q_{i_k}$ is non-empty, then $\Lambda(Q_{i_1}), \ldots , \Lambda (Q_{i_k})$ is a part of a basis for the integral lattice $\text{Hom}(S^1, T)\cong \Z^n$. Moreover, under our assumption of local standardness and $Q$ being a homology polytope, the manifold $X$ is determined up to equivariant diffeomorphisms by the pair $(Q, \Lambda)$ (see \cite[Lemma $4.5$]{masudapan}).

\begin{ex}\label{11}\rm
\begin{itemize}
\item[(1)] Let $\mathbb{T} \cong \left(\C^* \right)^n$ be the algebraic torus, $M=\text{Hom }(\mathbb{T}, \C^*)\cong  \text{Hom }(T, S^1)$ be the character lattice, and let $N=\text{Hom }(M, \Z)$ be the dual lattice. Consider the smooth complete $\mathbb{T}$-toric variety $X=X(\Delta)$ corresponding to a fan $\Delta$ in $N_{\R}:=N  \otimes_{\Z} \R \cong \R^n$ under the action of the torus $\mathbb{T}$. The orbit space of $X$ under the action of the compact torus $T \ (\subset \mathbb{T})$ is the manifold with corners $X_{\geq}$, which is formed by gluing $\left( U_{\sigma} \right)_{\geq}=\text{Hom}_{\text{sg}}(\sigma^{\vee} \cap M, \R_{\geq })$ $($see \cite[Section $4.1$]{ful}$)$. For each $\rho \in \Delta(1)$, let $v_{\rho} \in \text{Hom }(S^1, T)=N$ be the primitive ray generator of $\rho$. The characteristic submanifolds are given by the divisors ${\mathcal D}_{\rho}$ for $\rho \in \Delta(1)$, these are fixed by the circle subgroups $\text{Image}(v_{\rho})$. In this case the characteristic map $\Lambda$ is given by sending $({\mathcal D}_{\rho})_{\geq 0}$ to $v_{\rho}$. Since $H^{\ast}(X)$ is generated in degree two by \cite[Theorem $10.8$]{dan}, $X_{\geq}$ is a homology polytope. 
\item[(2)] Another class of examples are quasitoric manifolds introduced by Davis and \\ Januszkiewicz in \cite{davisjanu}. By definition a quasitoric manifold is locally standard under the $T$-action and the orbit space is a simple convex polytope and hence a homology polytope.
\end{itemize}
\end{ex}

\begin{rmk}{\rm
In \cite{suyama}, the author has constructed smooth complete toric varieties of complex dimension $\geq 4$ whose orbit spaces by the action of the compact torus are not homeomorphic to simple polytopes (as manifolds with corners). These provide the first known examples of smooth complete toric varieties that are not quasitoric manifolds.}
\end{rmk}

\noindent
The following lemma is an equivariant version of \cite[Lemma $5.1$]{param} and \cite[Proposition $2.1$]{uma}.
\begin{lemma}\label{1}
Let $X$ be a locally standard torus manifold with orbit space $X/T=Q$. For each $i$, $1 \leq i \leq d$, there exists a $T$-equivariant complex line bundle $L_i$ such that $c_1(L_i)=[V_i] \in H^2(X)$, where $[V_i]$ denotes the cohomology class dual to $V_i$ and each $L_i$ admits an equivariant section $s_i : X \rightarrow L_i$ which vanishes precisely along $V_i$. 
\end{lemma}
\noindent {\bf Proof:} Set $V=V_i$ and recall that $V$ is a closed $T$- invariant codimension $2$ submanifold of $X$ . Since $T$ is a compact Lie group we can assume that $X$ is endowed with a $T$-invariant Riemannian metric (see \cite[Chapter VI, Theorem $2.1$]{bredon}). Let $\nu$ denote the normal bundle to $V$ in $X$.  We have the decomposition $T(V) \oplus \nu=T(X)\mid_V$. Since $V$ is $T$-invariant, $T(V)$ and $T(X) \mid_V$ are $T$-equivariant vector bundles. Moreover, since the Riemannian metric is also $T$-invariant, $\nu=T(V)^{\perp} \subseteq T(X)\mid_V$ is naturally a $T$-equivariant real vector bundle. Furthermore, we see that $\nu$ is a canonically oriented real $2$-plane bundle since $T(V)$ and $T(X)\mid_{V}$ are oriented by the choice of the omniorientation. Thus $\nu$ admits a reduction of structure group to $SO(2,\mathbb{R})\cong S^1$ giving $\nu$ the structure of a complex line bundle. Since $T$ is a connected Lie group and $\nu$ is $T$-equivariant, $T$ preserves the orientation under the linear action on the fibre. (Fixing an oriented basis for the $\mathbb{R}$-vector space $\nu_x, ~ \mbox{for every} ~x\in V$, $t\mapsto \psi_t\in \text{Hom}(\nu_x, \nu_{tx})$ defines a continuous map from $T$ to $SO(2,\mathbb{R}) \subseteq O(2,\mathbb{R})$.) This implies that $T$ preserves the complex structure on the fibre, making $\nu$ a $T$-equivariant complex line bundle.
  
Now (by \cite[Chapter VI, Theorem $2.2$]{bredon}) $V$ has a closed invariant tubular neighbourhood denoted by $D$ which is equivariantly diffeomorphic to the disk bundle associated to the normal bundle $\nu$. The restriction of the equivariant diffeomorphism to the zero section of $\nu$ is the inclusion of $V$ in ${D}\subset X$. We denote by $\varpi : {D} \rightarrow V$ the projection map of the disk bundle. The complex line bundle $\varpi^*(\nu)$ admits an equivariant section $s: {D} \rightarrow \varpi^*(\nu)$ which vanishes precisely along $V$. Consider the trivial complex line bundle ${\mathcal E}:=(X \setminus \text{int}~ {D}) \times \C$ on $(X\setminus \text{int} ~{D})$, with the canonical $T$-action on $(X \setminus \text{int} ~{D})$ and the trivial $T$-action on the fibre $\C$. Consider the equivariant bundle isomorphism $\eta : {\mathcal E} \mid_{\partial {D}} \rightarrow \varpi^*(\nu)\mid_{\partial {D}}$ given by $(x, \lambda) \mapsto \lambda s(x)$, for all $x \in \partial {D}$. Now using clutching of bundles (see \cite[Theorem $3.2$]{karoubi}), glueing ${\mathcal E} \mid_{\partial {D}}$ along $\varpi^*(\nu)\mid_{\partial {D}}$ using the equivariant identification $\eta$ we get an equivariant line bundle, say $L$ on $X$. Note that $L$ admits an equivariant section $\tilde{s}$ (which restricts to $s$ on $D$ and $x \mapsto (x, 1)$ on $(X \setminus \text{int}~
  {D})$) that vanishes precisely along $V$. Hence $c_1(L)=[V]$ and this completes the proof.  $\hfill\square$

\begin{rmk}\label{12}\rm
  Let $p': ET \rightarrow BT$ be the universal principal $T$-bundle with the associated bundle $\pi' : ET \times_T X \rightarrow BT$. For a $T$-equivariant line bundle $q:L \rightarrow X$, we obtain the line bundle $ET \times_T L$ on $ET \times_T X$ with the projection $[e,l] \mapsto [e, q(l)]$. If $L$ has a $T$-invariant section $s$ which vanishes precisely along $V \subseteq X$, we obtain a section $\tilde{s}$ of $ET \times_T L$, defined by $[e,x] \mapsto [e, s(x)]$. Thus $\tilde{s}$ vanishes precisely along $ET \times_T V\subseteq ET\times_{T} X$. It follows that $c_1^T(L)=c_1(ET \times_T L)=[ET \times_T V]:=[V]_T$.  \end{rmk}

\begin{rmk}\label{hompol} {\em Note that in Lemma \ref{1} we do not assume that $Q$ is a homology polytope or even face-acyclic. It holds when $Q$ is simply a nice manifold with corners.} \end{rmk}

\begin{rmk} {\em Throughout this text by $H^*(~~)$ we shall always mean cohomology ring with $\mathbb{Z}$-coefficients unless specified otherwise.}\end{rmk}

\subsection{Cohomology ring and $K$-ring of torus manifolds} We now recall the presentation of the cohomology ring and $K$-ring of torus manifolds from \cite{masudapan} and \cite{param}.  \begin{thm}\label{cohom} $($\cite[Corollary $7.8$]{masudapan}, \cite[Proposition $5.2$]{param}$)$ Let $I$ be the ideal in $\Z[x_1, \ldots, x_d]$ generated by the elements: \begin{enumerate} \item[(i)] $x_{i_1} \cdots x_{i_r}$ whenever $V_{i_1} \cap \cdots \cap V_{i_r} = \emptyset$, \item[(ii)] $\displaystyle{\sum_{1 \leq i \leq d} \langle u, v_i \rangle x_i}$ where $u \in \text{Hom}(T, S^1)$.  \end{enumerate} We have an isomorphism of $\mathbb{Z}$-algebras \(\displaystyle\frac{\Z[x_1, \ldots, x_d]}{I}\stackrel{\sim}{\rightarrow} H^{\ast}(X)\) which maps $x_i$ to $c_1(L_i)=[V_i]\in H^2(X)$ for $1\leq i\leq d$. Furthermore, by \cite[Equation $(5.2)$, Section $7.2$]{masudapan}, $H^{\ast}(X)$ is a free abelian group of rank $\chi(X)=\mid X^T\mid=m$. Here $m$ equals the number of vertices of $Q$. \end{thm}

\begin{thm}\label{kring} \cite[Theorem $5.3$]{param} Let $J'$ be the ideal in $\Z[x_1, \ldots, x_d]$ generated by the following elements:
\begin{itemize}
\item[(i)] $x_{i_1} \cdots x_{i_r}$, whenever $V_{i_1} \cap \cdots \cap V_{i_r} = \emptyset$,
\item[(ii)]  $\displaystyle{\prod_{\{1 \leq i \leq d: \langle u, v_i \rangle > 0 \}} (1-  x_{i})^{ \langle u, v_i \rangle}- \prod_{\{1 \leq j \leq d: \langle u, v_j \rangle < 0 \}} (1-  x_{j})^{-  \langle u, v_j \rangle}}$ for $u \in \text{Hom}( T, S^1)$.
\end{itemize}
We have an isomorphism of $\mathbb{Z}$-algebras \(\displaystyle\frac{\Z[x_1, \ldots, x_d]}{J'} \stackrel{\sim}{\rightarrow} K^*(X)\) which maps $x_i$ to $1-[L_i]$, $1 \leq i \leq d$. Furthermore, $K^*(X)$ is a free abelian group of rank equal to $\chi(X)=m$ $($see \cite[Remark 3.2]{param}$)$.

\end{thm}

\begin{rmk}\label{10} \rm Let $J$ the ideal in $\mathbb{Z}[y^{\pm 1}_1,\ldots, y^{\pm 1}_d]$ generated by the following elements: \begin{itemize} \item[(i)] $\displaystyle{\prod_{1 \leq j \leq r} \left(1- y_{i_j} \right)}$, whenever $V_{i_1} \cap \cdots \cap V_{i_r} = \emptyset$, \item[(ii)] $\displaystyle{\prod_{1 \leq i \leq d} y_i^{\langle u, v_{i} \rangle}}$ for $u \in \text{Hom}( T, S^1)$.  \end{itemize} In Theorem \ref{kring}, by making the transformation $y_i=1-x_i$, $1 \leq i \leq d$ we get the following alternative presentation \(\displaystyle \frac{\Z[y_1^{\pm 1}, \ldots, y_d^{\pm 1}]}{J}\) for $K^*(X)$ which sends $y_i$ to $[L_i]$, $1 \leq i \leq d$ (see \cite[Remark 4.2]{param}). \end{rmk}

\section{Cohomology ring of torus manifold bundles}

Let $p: E\rightarrow B$ be a principal bundle with fibre and structure group the compact torus $T$ over a topological space $B$. Then one has the associated fibre bundle $\pi :  E(X)\rightarrow B$ with fibre the torus manifold $X$, where $E(X):= E \times_T X $, and $\pi ([e,x])=p(e)$. For $u \in \text{Hom }(T,S^1)$, let $\C_u$ denote the corresponding $1$-dimensional $T$-representation. One has a $T$-line bundle $\xi_u$ on $B$ whose total space is $E \times_T \C_u$. For the $T$-equivariant complex line bundle $L_i$ from Lemma \ref{1}, let $E(L_i):=E \times_T L_i$ denote the associated line bundle on $E(X)$.

\begin{defn}\label{3}\rm Let $R(B, \left( Q, \Lambda \right))$ denote the ring \( \displaystyle \frac{H^*(B)[x_1, \ldots, x_d]}{\mathcal{I}} \) where the ideal $\mathcal{I}$ is generated by the following elements:
\begin{itemize}
\item[(i)] $x_{i_1} \cdots x_{i_r}$, whenever $Q_{i_1} \cap \cdots \cap Q_{i_r} = \emptyset$,
\item[(ii)] $\displaystyle{\sum_{i=1 }^d \langle u, v_i \rangle x_{i}-c_1(\xi_u)}$ for $u \in \text{Hom}( T, S^1)$.
\end{itemize}
\end{defn}

\noindent
Recall that the \emph{face ring} $\Z[Q]$ of the homology polytope $Q$ is defined to be \( \displaystyle \frac{\Z[x_1, \ldots, x_d]}{I_1} \) where the ideal $I_1$ is generated by elements of the form 
\begin{equation}
x_{i_1} \cdots x_{i_r}, \text{ whenever } Q_{i_1} \cap \cdots \cap Q_{i_r} = \emptyset.
\end{equation}

For $h(x_1,\ldots,x_d)\in \mathbb{Z}[x_1,\ldots,x_d]\subseteq H^*(B)[x_1,\ldots,x_d]$ we shall denote by $\bar{h}(x_1,\ldots,x_d)$ its class in $\mathbb{Z}[Q]$ and by $\bar{\bar{h}}(x_1,\ldots,x_d)$ its class in $R(B, \left( Q, \Lambda \right))$.

We have a canonical $H^{\ast}(BT)$-algebra structure on $\Z[Q]$ given by the ring homomorphism $H^{\ast}(BT) \rightarrow \Z[Q]$ which maps $u$ to $\displaystyle{\sum_{i=1}^d \langle u, v_i \rangle \bar{x_i}}$. A canonical $H^*(B)$-module structure on $H^*(B) \otimes_{H^{\ast}(BT)} \Z[Q]$ is obtained by extending scalars to $H^*(B)$ via the homomorphism $H^*(BT) \rightarrow H^*(B)$ that sends $u$ to $c_1(\xi_u)$.

\begin{lemma}\label{4} We have an isomorphism $R(B, \left( Q, \Lambda \right)) \cong H^*(B) \otimes_{H^*(BT)} \Z[Q]$ of $H^*(B)$-modules. In particular, $R(B, \left( Q, \Lambda \right))$ is a free $H^*(B)$-module of rank $m$.  \end{lemma}

\noindent
{\bf Proof:} Define $\alpha : H^*(B)[x_1,\ldots,x_d] \rightarrow  H^*(B) \otimes_{H^*(BT)} \Z[Q]$ by sending $x_i \mapsto 1 \otimes \bar{x_i}$ and $b \mapsto b \otimes 1$ for $b \in H^*(B)$. Clearly the generators of $\mathcal{I}$ listed in $(i)$ of Definition \ref{3} map to zero under $\alpha$. Now \[\alpha(\sum_{i=1 }^d \langle u, v_i \rangle x_{i}-c_1(\xi_u))=1 \otimes \sum_{i=1 }^d \langle u, v_i \rangle \bar{x_{i}}- c_1(\xi_u) \otimes 1=1 \otimes u \cdot 1- u \cdot 1 \otimes 1\] which is zero in $ H^*(B) \otimes_{H^*(BT)} \Z[Q]$. Hence $\alpha$ induces a well defined $H^*(B)$-module homomorphism $\bar{\alpha} : R(B, \left( Q, \Lambda \right))\rightarrow H^*(B) \otimes_{H^*(BT)} \Z[Q]$.\\
We define $\beta : H^*(B) \otimes_{\Z}\Z[Q] \rightarrow R(B, \left( Q, \Lambda \right))$ by $b \otimes \bar{h}(x_1, \ldots, x_d) \mapsto b \bar{\bar{h}}(x_1, \ldots, x_d)$, for $b \in H^*(B)$ and $\bar{h}(x_1, \ldots, x_d) \in \Z[Q]$. Clearly this is well defined. Now \[\beta(1 \otimes u \cdot 1- u \cdot 1 \otimes 1)=\sum_{i=1 }^d \langle u, v_i \rangle \bar{\bar{x_i}}-c_1(\xi_u)\] which is zero in $R(B, \left( Q, \Lambda \right))$. Hence $\beta$ induces a map $\bar{\beta} : H^*(B) \otimes_{H^*(BT)} \Z[Q] \rightarrow R(B, \left( Q, \Lambda \right))$. Noting that $\bar{\alpha}$ and $\bar{\beta}$ are inverses of each other proves the first assertion. Now $\Z[Q]$ is free $H^*(BT)$-module of rank $m$ by \cite[Theorem $7.7$, Lemma $2.1$]{masudapan}, which proves the second assertion. $\hfill\square$

\noindent The following is the main theorem of this section.  \begin{thm}\label{2} Let $B$ have the homotopy type of a finite CW complex. The map $\Phi : R(B, \left( Q, \Lambda \right)) \rightarrow H^*(E(X))$ which sends $x_i$ to $c_1(E(L_i))$ is an isomorphism of $H^*(B)$-algebras.
\end{thm}

\noindent
{\bf Proof: }Suppose that $Q_{i_1} \cap \cdots \cap Q_{i_r} = \emptyset$, which implies $V_{i_1} \cap \cdots \cap V_{i_r} = \emptyset$. So by Lemma \ref{1}, the bundle $L_{i_1}\oplus \cdots \oplus L_{i_r}$ has a nowhere vanishing $T$-equivariant section. Hence by Remark \ref{12}, the bundle $E(L_{i_1})\oplus \cdots \oplus E(L_{i_r})$ admits a nowhere vanishing section. This shows that 
\begin{equation}
c_1(E(L_{i_1}) \cdots c_1( E(L_{i_r}))=0
\end{equation}
in $H^{2r}(E(X))$. Hence the elements listed in (i) of Definition \ref{3} map to zero under $\Phi$.

Let $L_u$ be the trivial line bundle $X \times \C_u$ on $X$. Consider the associated line bundle ${\xi}'_u:=ET \times_T \C_u$ on $BT$. Note that $ET \times_T L_u$ is isomorphic to the pullback $\pi'^*({\xi}'_u)$ where $\pi'$ is as in Remark \ref{12}. By the naturality of Chern classes \begin{equation}\label{equivchern} c_1^{T}(L_u):=c_1(ET \times_T L_u)=\pi'^*(c_1({\xi}'_u)).\end{equation} By \cite[Proposition $3.3$]{masudapan},  \begin{equation}\label{equivchern1}\pi'^*(c_1({\xi}'_u))=\sum_{i=1}^d \langle u, v_i \rangle [V_i]_T\end{equation} and by Remark \ref{12},  \begin{equation}\label{equivchern2} c^T_1\left( \prod_{i=1}^d L_i^{\langle u, v_i \rangle}\right)=\sum_{i=1}^d \langle u, v_i \rangle [V_i]_T\end{equation} in $H^2_{T}(X)$. Now, (\ref{equivchern}), (\ref{equivchern1}) and (\ref{equivchern2}) together imply $\displaystyle {c_1^{T}(\prod_{i=1}^d L_i^{\langle u, v_i \rangle}) =c_1^{T}(L_u)}$.
This in turn implies that $\displaystyle {L_u \cong \prod_{i=1}^d L_i^{\langle u, v_i \rangle}}$ as $T$-equivariant line bundles by \cite[Theorem C$.47$]{moment}. Thus  \begin{equation}\label{associatediso}\pi^*(\xi_u) \cong E(L_u) \cong \prod_{i=1}^d E(L_i)^{\langle u, v_i \rangle} .\end{equation} Taking first Chern classes on both sides of (\ref{associatediso}) we get 
\begin{equation}\label{chernassociatediso}
\sum_{i=1}^d \langle u, v_i \rangle c_1(E(L_i))=c_1(\pi^*(\xi_u)).
\end{equation}
This implies that the generators of $\mathcal{I}$ listed in (ii) of Definition \ref{3} map to zero under $\Phi$, hence it is a well-defined ring homomorphism.

By Theorem \ref{cohom}, there exist $p_i(x_1, \ldots, x_d) \in \Z[x_1, \ldots, x_d]$, $1 \leq i \leq m$ such that \[p_i:=p_i(c_1(L_1),\ldots, c_1(L_d)) :1\leq i\leq m\] form a $\Z$-basis of $H^*(X)$. Consider \[ \displaystyle {P}_i:=p_i(c_1(E(L_1)),\ldots, c_1(E(L_d))) :1\leq i\leq m\] in $H^*(E(X))$. Since $E(L_i)|_X=L_i$, it follows that ${P}_j\mid_X=p_j$. Since $H^k(X)$ is free for all $k$, by the Leray-Hirsch theorem, (see \cite[Theorem $4D.1)$]{top}) $H^*(E(X))$ is a free $H^*(B)$-module with ${P}_1, \ldots, {P}_m$ as a basis. Moreover, since $\Phi(x_i)=c_1(E(L_i))$, each ${P}_i$ has a preimage under $\Phi$. Hence by Lemma \ref{4}, $\Phi$ is a surjective $H^*(B)$-module map between two free $H^*(B)$-modules of the same rank.  Furthermore, since $H^*(B)$ is a finitely generated abelian group, it follows that $\Phi$ is a surjective map from a finitely generated abelian group to itself, and hence an isomorphism. (More generally, a surjective morphism from a finitely generated module over a Noetherian commutative ring to itself is an isomorphism (see \cite[Chapter $6$, Exercise $1.(i)$]{AM})).  $\hfill\square$

\section{$K$-ring of torus manifold bundles}

\begin{defn}\label{14} \rm  The $K$-theoretic face ring of the homology polytope $Q$ is defined to be \( \displaystyle {\mathcal K}(Q):=\frac {\Z[y_1^{\pm 1}, \ldots, y_d^{\pm 1}]}{J_1} \) where $J_1$ is the ideal generated by elements of the form 
\begin{equation}
\left(1- y_{i_1} \right) \cdots \left(1- y_{i_r} \right), \text{ whenever } Q_{i_1} \cap \cdots \cap Q_{i_r} = \emptyset.
\end{equation}
 \end{defn}

\noindent
We show that ${\mathcal K}(Q)$ is a free $RT$-module in Proposition \ref{8}. We first set up the notation. Recall that $v_i \in \text{Hom }(S^1, T)$ determines the circle subgroup of $T$ fixing $V_i$ for $i=1, \ldots,d$. Let $\mathcal{V}$ denote the set of vertices of $Q$ and let $a \in \mathcal{V}$. Write $a=Q_{i_1} \cap \cdots \cap Q_{i_n}$ as an intersection of facets. Then by \cite[Proposition $3.3$]{masudapan}, the elements $v_{i_1}, \ldots, v_{i_n}$ form a basis of $\text{Hom }(S^1,T)$.  We set \[RT_a:=\Z[\chi^{\pm u_{i_1}}, \ldots, \chi^{\pm u_{i_n}}]\] where $u_{i_1}, \ldots, u_{i_n}$ denotes the dual basis of $v_{i_1}, \ldots, v_{i_n}$. For any $b \in \mathcal{V}$, denote by $a \vee b$ the minimal face of $Q$ containing both $a$ and $b$. If $a \vee b=Q$, set $RT_{a \vee b}=\Z$ and the projection map $RT_a \rightarrow RT_{a \vee b}$ to be the augmentation map. Otherwise when $a \vee b$ is a proper face, write $a \vee b=Q_{i_1} \cap \cdots \cap Q_{i_l}$ and set \[ RT_{a \vee b}:=\Z[\frac{M}{\langle v_{i_1}, \ldots, v_{i_l} \rangle^{\perp}}]=\Z[\chi^{\pm u_{i_1}}, \ldots, \chi^{\pm u_{i_l}}]. \] Then we have the canonical projection map $RT_a \rightarrow RT_{a \vee b}$ given by $\chi^{u_{i_j}} \mapsto \chi^{u_{i_j}}$ for $j=1, \ldots, l$ and $\chi^{u_{i_j}} \mapsto 1$ for $j= l+1, \ldots, n$.

The following lemma is analogous to \cite[Theorem $6.4$]{vezzosi2003higher} in the setting of torus manifolds. We prove it along similar lines.

\begin{lemma}\label{7}
There is an inclusion of rings \[\displaystyle \bar{\phi}: {\mathcal K}(Q) \hookrightarrow \prod_{a \in \mathcal{V}} RT_a. \] The image consists of elements of the form $\displaystyle{\left(r_a \right) \in \prod_{a \in \mathcal{V}} RT_a}$, where for any two distinct $a, \ b \in \mathcal{V}$, the restriction of $r_a$ and $r_b$ to $RT_{a \vee b}$ coincide.
\end{lemma}
\noindent
{\bf Proof: } Define the map $\displaystyle{\phi:\Z[y_{1}^{\pm 1}, \ldots, y_{d}^{\pm 1}] \rightarrow  \prod_{a \in \mathcal{V}} RT_a}$ given by $y_i \mapsto r_i:=(r_{ia})$ where 
\[ r_{ia} = \left\{ \begin{array}
{r@{\quad \quad}l}
1 & \text{if } a \notin Q_i \\ 
\chi^{u_i} &  \text{if } a \in Q_i \\
\end{array} \right. \]
Set $\displaystyle {W=\{\left(r_a \right)  \in \prod_{a \in \mathcal{V}} RT_a: r_a \mid_{a \vee b}=r_b \mid_{a \vee b}, \text{for all } a \neq b \in \mathcal{V} \}}$. Note that $W$ is a subring of $\displaystyle{\prod_{a \in \mathcal{V}} RT_a}$.

Let $a, b \in \mathcal{V}$ be distinct. If $a \vee b=Q$, then there is nothing to prove. Otherwise write $a \vee b =Q_{i_1} \cap \cdots \cap Q_{i_l}$, where $a=Q_{i_1} \cap \cdots \cap Q_{i_n}$ and $b=Q_{i_1} \cap \cdots \cap Q_{i_l} \cap Q_{j_{l+1}} \cap \cdots \cap Q_{j_n}$. Now consider the following cases:
\begin{enumerate}
\item $a, b \notin Q_i$: Then $r_{ia}=1=r_{ib}$, hence $r_a \mid_{a \vee b}=r_b \mid_{a \vee b}$.
\item $a \notin Q_i$ and $b \in Q_i$: Then $r_{ia}=1$ and $r_{ib}=\chi^{u_i}$. Note that under the restriction map $RT_b \rightarrow RT_{a \vee b}$, $\chi^{u_i} \mapsto 1$, since $u_i \in \langle v_{i_1}, \ldots, v_{i_l} \rangle^{\perp}$. Hence we are done in this case.
\item $a, b \in Q_i$: Then $r_{ia} \mid_{a \vee b}=\chi^{u_i}=r_{ib} \mid_{a \vee b}$ and under the respective projection they map to the same image since $i \in \{i_1, \ldots, i_l \}$.
\end{enumerate} 
This proves that $r_i \in W$ for $1 \leq i \leq d$. We show that elements of $W$ can be written as Laurent polynomials in $r_i$'s.
Set \[\mathcal{V}=\{a_1, \ldots, a_m \}\] and let $\alpha=\left( \alpha_{a_i} \right) \in W$. Let $a_1=Q_{i_1} \cap \cdots \cap Q_{i_n}$, then $\alpha_{a_1} \in RT_{a_1}=\Z[\chi^{\pm u_{i_1}}, \ldots, \chi^{\pm u_{i_n}}]$ and hence we can find a Laurent polynomial $p_1(y_{i_1}, \ldots, y_{i_n})$ such that $p_1(r_{i_1}, \ldots, r_{i_n})_{a_1}=\alpha_{a_1}$. Let $\alpha_1:= \alpha-p_1(r_{i_1}, \ldots, r_{i_n})$. Then we see that $\alpha_{1_{a_1}}=0$.

Now let $a_2=Q_{i_1} \cap \cdots \cap Q_{i_l} \cap Q_{j_{l+1}} \cap \cdots \cap Q_{j_n}$ such that $a_1 \vee a_2= Q_{i_1} \cap \cdots \cap Q_{i_l}$. Similarly as above there is a Laurent polynomial $p_2(y_{i_1}, \ldots, y_{i_l}, y_{j_{l+1}}, \ldots, y_{j_n})$ such that $p_2(r_{i_1}, \ldots, r_{i_l}, r_{j_{l+1}}, \ldots, r_{j_n})_{a_2}=\alpha_{1_{a_2}}$. Note that \[p_2(r_{i_1}, \ldots, r_{i_l}, r_{j_{l+1}}, \ldots, r_{j_n})_{a_1}=p_2(\chi^{u_{i_1}}, \ldots, \chi^{u_{i_l}}, 1, \ldots, 1)\]
whose projection to $RT_{a_1 \vee a_2}$ remains unchanged, i.e. \begin{equation}\label{eq1} p_2(r_{i_1}, \ldots, r_{i_l}, r_{j_{l+1}}, \ldots, r_{j_n})_{a_1}= p_2(r_{i_1}, \ldots, r_{i_l}, r_{j_{l+1}}, \ldots, r_{j_n})_{a_1} \mid _{a_1 \vee a_2}.\end{equation}

Since $\alpha_1\in W$, $\alpha_{1_{a_2}} \mid_{a_1 \vee a_2}=\alpha_{1_{a_1}} \mid_{a_1 \vee a_2}=0$. Moreover, $p_2(r_{i_1}, \ldots, r_{i_l}, r_{j_{l+1}}, \ldots, r_{j_n}) \in W$ implies \[ p_2(r_{i_1}, \ldots, r_{i_l}, r_{j_{l+1}}, \ldots, r_{j_n})_{a_1} \mid _{a_1 \vee a_2}= p_2(r_{i_1}, \ldots, r_{i_l}, r_{j_{l+1}}, \ldots, r_{j_n})_{a_2} \mid _{a_1 \vee a_2}=\alpha_{1_{a_2}}\mid_{a_1\vee a_2}=0.\] Now, (\ref{eq1}) implies \[p_2(r_{i_1}, \ldots, r_{i_l}, r_{j_{l+1}}, \ldots, r_{j_n})_{a_1}=0.\] Letting $\alpha_2 := \alpha_1-p_2(r_{i_1}, \ldots, r_{i_l}, r_{j_{l+1}}, \ldots, r_{j_n})$, we have $\alpha_{2_{a_1}}=0=\alpha_{2_{a_2}}$.  Repeating this process for $a_3, \ldots, a_m$, where $a_k=Q_{k_1}\cap\cdots\cap Q_{k_n}$ for $k=1,\ldots, m$, we get that $\alpha_{m_{a_1}}=\alpha_{m_{a_2}}=\cdots=\alpha_{m_{a_m}}=0$, for $\alpha_m=\alpha-\sum_{k=1}^m p_k(r_{k_1},\ldots,r_{k_n})$. Thus $\alpha_m=0$, so that $\alpha$ is in the image of $\phi$. Since $\alpha\in W$ was arbitrary, $\phi$ is surjective. It remains to show that $\text{ker}(\phi)=J_1$.

For $a \in \mathcal{V}$, consider the map $\phi_a:\Z[y_{1}^{\pm 1}, \ldots, y_{d}^{\pm 1}] \rightarrow RT_a$ which sends $y_i \mapsto r_{ia}$ for $1 \leq i \leq d$. We see that $\text{ker}(\phi_a)=J_a:=\langle y_j-1: a \notin Q_j \rangle $ and clearly $\text{ker}(\phi)= \cap_{a \in \mathcal{V}} J_a$. Then $\cap_{a \in \mathcal{V}} J_a= J_1$ follows from \cite[Lemma $6.5$]{vezzosi2003higher}. Hence we get the induced ring homomorphism $\displaystyle{\bar{\phi}:\mathcal{K}(Q)\stackrel{\sim}{\rightarrow} W\hookrightarrow \prod_{a\in \mathcal{V}} RT_a}$ as required.  $\hfill\square$

\vspace{1cm}

Note that one has a monomorphism of rings $RT \stackrel{\iota}{\rightarrow }{\mathcal K}(Q)$ defined by $\displaystyle{\chi^u \mapsto \prod_{1 \leq i \leq d} y_{i}^{\langle u, v_{i} \rangle}}$, $u \in \text{Hom}( T, S^1)=M$, which gives an $RT$-algebra structure on ${\mathcal K}(Q)$.

Moreover, for every $a_k=Q_{k_1}\cap\cdots \cap Q_{k_n}$, in $\mathcal{V}$, we have the isomorphism $\zeta_{k}:\mathbb{Z}[M]=RT \rightarrow RT_{a_k}$ which maps $\displaystyle{\chi^u\mapsto \prod_{j=1}^n\chi^{\langle u,v_{k_j}\rangle u_{k_j}}}$ for $1\leq k\leq m$. Thus $\displaystyle{\prod_{k=1}^m \zeta_k}$ identifies $(RT)^m$ with $\displaystyle{\prod_{a\in \mathcal{V}} RT_a=\prod_{k=1}^m RT_{a_k}}$. Now, $(RT)^m$ has a canonical $RT$-algebra structure via the diagonal embedding $\delta$. Hence $\displaystyle{\zeta=(\prod_{k=1}^m\zeta_k)\circ \delta:RT\longrightarrow \prod_{k=1}^m RT_{a_k}}$ which maps $\displaystyle{\chi^u\mapsto (\prod_{j=1}^n\chi^{\langle u,v_{k_j}\rangle u_{k_j}})}$ gives the canonical $RT$-algebra structure on $\displaystyle{\prod_{a\in \mathcal{V}} RT_a}$.

\begin{cor}\label{algegramono}
The inclusion of rings $\bar{\phi}$ in {\em Lemma \ref{7}} is a monomorphism of $RT$-algebras.
\end{cor}
{\bf Proof:} The proof follows readily since it can be seen that $\bar{\phi}\circ \iota=\zeta$.$\hfill\square$

\begin{prop}\label{8} 
${\mathcal K}(Q)$ is a free $RT$-module of rank $\chi(X)$.
\end{prop}

\noindent
{\bf Proof: } We see that ${\mathcal K}(Q)$  is isomorphic to a localization of $\Z[Q]$ by a similar argument as in the proof of \cite[Theorem $2.3$]{baggio2007equivariant}. Explicitly, we see that there is an ring isomorphism 
\( \displaystyle \Z[Q] \cong \frac{\Z[y_1, \ldots, y_d]}{J_1 \cap \Z[y_1, \ldots, y_d] }\)
which sends $x_i$ to $y_i-1$. This remains an isomorphism if we localize at the respective multiplicative systems $S_I=\{(x_i+1)^k\}_{k \in \N}$ and $S_J=\{y_i^k\}_{k \in \N}$:
\[S_I^{-1}\Z[Q] \cong S_J^{-1}\frac{\Z[y_1, \ldots, y_d]}{J_1 \cap \Z[y_1, \ldots, y_d] }={\mathcal K}(Q)\] Now $\Z[Q]$ is Cohen-Macaulay by \cite[Lemma $8.2$]{masudapan}, hence ${\mathcal K}(Q)$ is also Cohen-Macaulay. Note that ${\mathcal K}(Q)$ is a finite $RT$-module since it is a submodule of a Noetherian module $\prod_{a \in \mathcal{V}} RT_a\simeq RT^m$ by Lemma \ref{7}.  Hence $\iota:RT \subseteq {\mathcal K}(Q)$ is an integral extension. Since $RT$ is an integrally closed domain and ${\mathcal K}(Q)$ is a torsion free $RT$-module, by the Going Down Theorem (\cite[Corollary $2.2.8$]{huneke}), for any maximal ideal $\mathfrak{M}$ of ${\mathcal K}(Q)$ which contracts to the maximal ideal $\mathfrak{m}$ of $RT$, $\text{ht }{\mathfrak{m}}= \text{ht }{\mathfrak{M}}$. Then by \cite[Lemma $2.4$]{baggio2007equivariant}, ${\mathcal K}(Q)$ is a projective $RT$-module. Moreover, since $RT$ is a Laurent polynimial ring, ${\mathcal K}(Q)$ is in fact a free $RT$-module. Now note that the presentation of $K^*(X)$ in \cite[Theorem 5.3]{param} and Remark \ref{10}, implies that $$K^*(X) \cong \frac {\Z[y_1^{\pm 1}, \ldots, y_d^{\pm 1}]}{J} \cong \mathbb{Z}\otimes_{RT} {\mathcal K}(Q) $$ where the extension of scalars to $\mathbb{Z}$ is via the augmentation homomorphism $RT \stackrel{\epsilon}{\rightarrow} \mathbb{Z}$. On the other hand it is also known that $K^*(X)$ is a free abelian group of rank $\chi(X)$. Hence the proposition follows.  $\hfill\square$

\begin{defn}\label{5}\rm Let \(\displaystyle\mathcal{R}(B, \left( Q, \Lambda \right)):= \frac{K^*(B)[y_1^{\pm 1}, \ldots, y_d^{\pm 1}]}{\mathcal{J}} \) where the ideal $\mathcal{J}$ is generated by the following elements:
\begin{itemize}
\item[(i)] $(1-y_{i_1}) \cdots (1-y_{i_r})$, whenever $Q_{i_1} \cap \cdots \cap Q_{i_r} = \emptyset$,
\item[(ii)] $\displaystyle{\prod_{1 \leq i \leq d} y_i^{\langle u, v_{i} \rangle}-[\xi_u]}$ for $u \in \text{Hom}( T, S^1)$.
\end{itemize}

\end{defn}

Consider the ring $K^*(B) \otimes_{RT} {\mathcal K}(Q)$ obtained from the  $RT$-algebra $\mathcal{K}(Q)$ by extending scalars to $K^*(B)$ via the homomorphism $RT \rightarrow K^*(B)$ which maps $\chi^u \mapsto [\xi_u]$. In particular,  by Proposition \ref{8}, $K^*(B) \otimes_{RT} {\mathcal K}(Q)$ is a free $K^*(B)$-module of rank $\chi(X)=m$.

\begin{lemma}\label{9} We have an isomorphism
\( \displaystyle \mathcal{R}(B, \left( Q, \Lambda \right)) \cong K^*(B) \otimes_{RT} {\mathcal K}(Q)\) as $K^*(B)$-modules. In particular, $\mathcal{R}(B, \left( Q, \Lambda \right))$ is a free $K^{\ast}(B)$-module of rank $\chi(X)$. 
\end{lemma} 
 
The proof is similar to the proof of Lemma \ref{4}.

\begin{thm}\label{6} Let $B$ has the homotopy type of a finite CW complex. Then we have an isomorphism \( \displaystyle \Psi:\mathcal{R}(B, \left( Q, \Lambda \right)) \stackrel{\sim}{\rightarrow} K^*(E(X)) \) of $K^*(B)$-modules, which maps $ y_i \mapsto [E(L_i)]$.
\end{thm}

Suppose that $Q_{i_1} \cap \cdots \cap Q_{i_r} = \emptyset$. Recall from the proof of the Theorem \ref{2}, the bundle $E(L_{i_1})\oplus \cdots \oplus E(L_{i_r})$ admits a nowhere vanishing section. Then applying $\gamma^r$-operation, we obtain $\gamma^r([L_{i_1}\oplus \cdots \oplus L_{i_r}]-r)=(-1)^r c_r(L_{i_1}\oplus \cdots \oplus L_{i_r})=0$. Also note that $\displaystyle{\gamma^r([L_{i_1}\oplus \cdots \oplus L_{i_r}]-r)=\prod_{1 \leq j \leq r} \left([L_{i_j}]-1 \right)}$. This shows that the elements listed in (i) of Definition \ref{5} map to zero under $\Psi$.

Note that, we have $\displaystyle{\pi^*(\xi_u) \cong E(L_u) \cong \prod_{i=1}^d E(L_i)^{\langle u, v_i \rangle}}$ from the proof of Theorem \ref{2}. This implies that the generators of $\mathcal{J}$ listed in (ii) of Definition \ref{5} map to zero under $\Psi$.

The surjectivity of $\Psi$ follows from the same argument as in the proof of Theorem \ref{2}, using a version of the Leray-Hirsch theorem in the setting of K-theory (see \cite[Theorem $2.25$]{hatcher2003vector}).

Then using Lemma \ref{9}, similar arguments as in the proof of Theorem \ref{2} show that $\Psi$ is an isomorphism.
$\hfill\square$

\section{Some applications}

As an illustration of the above results, we derive both the cohomology and $K$-ring of $E(X)$, where $X=X(\Delta)$ is a smooth complete toric variety (see Example \ref{11}).

\begin{defn}\rm For a smooth complete fan $\Delta$ we define the following rings.
\begin{enumerate}

\item Let $R(H^*(B), \Delta)$ denote the ring \( \displaystyle \frac{H^*(B)[X_1, \ldots, X_d]}{\mathcal{I}} \) where the ideal $\mathcal{I}$ is generated by the following elements:
\begin{itemize}
\item[(i)] $X_{i_1} \cdots X_{i_r}$, whenever $\rho_{i_1}, \ldots, \rho_{i_r}$ do not generate a cone in $\Delta$,
\item[(ii)] $\displaystyle{\sum_{i=1 }^d \langle u, v_i \rangle X_{i}-c_1(\xi_u)}$ for $u \in \text{Hom}( T, S^1)$.
\end{itemize}

\item Let $\mathcal{R}(K^*(B), \Delta)$ denote the ring \( \displaystyle \frac{K^*(B)[Y_1^{\pm 1}, \ldots, Y_d^{\pm 1}]}{\mathcal{J}} \) where the ideal $\mathcal{J}$ is generated by the following elements:
\begin{itemize}
\item[(i)] $(1-Y_{i_1}) \cdots (1-Y_{i_r})$, whenever $\rho_{i_1}, \ldots, \rho_{i_r}$ do not generate a cone in $\Delta$,
\item[(ii)] $\displaystyle{\prod_{1 \leq i \leq d} Y_i^{\langle u, v_{i} \rangle}-[\xi_u]}$ for $u \in \text{Hom}( T, S^1)$.
\end{itemize}

\end{enumerate}
\end{defn}

\begin{cor}\label{15} Let $X=X(\Delta)$ be a smooth complete $\mathbb{T} \cong (\C^*)^n$-toric variety. Let $p: E \rightarrow B$ be a principal $\mathbb{T}$-bundle, where $B$ has the homotopy type of a finite CW complex. 

\begin{enumerate}
\item[$(1)$]  The cohomology ring of $E(X)$ is isomorphic as an $H^*(B)$-algebra to $R(H^*(B), \Delta)$ under the isomorphism \( \displaystyle \Phi : R(H^*(B), \Delta) \rightarrow H^*(E(X))\) which sends $X_i$ to $c_1(E(L_i))$.

\item[$(2)$] The topological $K$-ring of $E(X)$ is isomorphic as a $K^*(B)$-algebra to $\mathcal{R}(K^*(B), \Delta)$ under the isomorphism \( \displaystyle \Psi : \mathcal{R}(K^*(B), \Delta)\rightarrow K^*(E(X))\) which sends $Y_i$ to $[E(L_i)]$.

\end{enumerate}
\end{cor} {\bf Proof:} We consider $X$ as a torus manifold with locally standard $T \cong (S^1)^n$ action and orbit space the homology polytope $X_{\geq}$ (see Example \ref{11}). We then have the principal $T$-bundle $p' : E \rightarrow E/T$ since $T$ is an admissible subgroup of $\mathbb{T}$ (i.e. $\mathbb{T} \rightarrow \mathbb{T}/T$ is a principal $T$-bundle). Note that $E \times _T X$ and $E \times _{\mathbb{T}} X$ are homotopy equivalent, since $\mathbb{T}=T \times (\R_{\geq})^n$ and $(\R_{\geq})^n$ is contractible. Similarly $B$ and $E/T$ are homotopy equivalent. The assertions $(1)$ and $(2)$ of the corollary now follow by applying Theorem \ref{2} and Theorem \ref{6} respectively for the space $E\times_{T} X$ associated to the principal $T$-bundle $p' : E \rightarrow E/T$.  Note that in the proof of assertion $(2)$ above, Proposition \ref{8} is immediate from \cite[Theorem $6.9$]{vezzosi2003higher} because the ring ${\mathcal K}(Q)$ in Proposition \ref{8} is the algebraic $\mathbb{T}$-equivariant $K$-ring of $X$. $\hfill\square$

\begin{rmk}{\rm Let $X$ be a torus manifold with locally standard action and orbit space a homology polytope $Q$ whose nerve is a shellable simplicial complex (see \cite{uma}), e.g. quasitoric manifolds. Then Theorem \ref{2} (respectively, Theorem \ref{6}) can be proved for for $B$ any topological space (respectively, $B$ compact Hausdorff topological space) using \cite[Lemma $2.1$, Lemma $2.2$]{paramuma}.  In particular, this gives a relative version of \cite[Theorem 1.3]{uma} and a generalization of \cite[Theorem 1.2]{MR2313027}.} \end{rmk}

\section{Torus manifold bundles when $X/T$ is not a homology polytope}

In the preceeding sections we considered torus manifolds $X$ defined in Section 2.1 with the additional assumption that $X/T=Q$ is a homology polytope. This ensured that the cohomology ring $H^*(X)$ was generated by the degree $2$ classes corresponding to the fundamental classes $[V_i]$ of the characteristic submanifolds (see Theorem \ref{cohom}). 

In this section we shall consider a torus manifold $X$ with a locally standard action of $T$ as defined in Section 2.1, with the exception that $X/T=Q$ is not assumed to be a homology polytope but only face acyclic. In particular, we do not assume that the prefaces are connected. Since $Q$ is face acyclic the cohomology ring of $X$ satisfies the property that $H^{odd}(X)=0$ (see \cite[Theorem 7.4.46]{bp}), which in particular also implies by the universal coefficient theorem that $H^*(X)$ is torsion free and hence free of rank $\chi(X)$. Moreover, \cite[Corollary 7.8]{masudapan} gives an explicit presentation of the ring $H^*(X)$. 

Let $p : E\longrightarrow B$ be a principal $T$-bundle and let $E(X):=E\times_{T} X$ be the associated torus manifold bundle. Let $B$ be a topological space having the homotopy type of a finite CW complex.

We then have the following theorem which gives a presentation of $H^*(E(X))$ as a $H^*(B)$-algebra.

\begin{thm}\label{cohombuntorus} Let $\mathfrak{I}$ be the ideal in the ring ${\mathfrak R}:= H^*(B) [x_{F}:F ~\mbox{a face of}~ Q]$ generated by the following relations:

(i) $\displaystyle{x_Gx_{H}-x_{G\vee H}\sum_{E\in G\cap H} x_{E}};$  
  
(ii) $\displaystyle{\sum_{i=1}^d \langle u,v_i\rangle x_{Q_i}}-c_1(\xi_u)$ for $u\in Hom(T,S^1)$ where $Q_i$ are the facets of $Q$, $v_i=\Lambda(Q_i)$ is the primitive vector in $Hom(S^1,T)\simeq \mathbb{Z}^n$ which determines the circle subgroup of $T$ fixing the characteristic submanifold $V_i$ for $1\leq i\leq d$, and $\xi_u=E\times_{T} \mathbb{C}_u$ is the line bundle on $B$ associated to the character $u\in M$. Since $X$ is omnioriented $v_i$ is well defined. {\em (See Section 2 and Section 3)}.
 
The map $\Phi_1 : {\mathfrak R} \rightarrow H^*(E(X))$ which sends $x_F$ to $[E(V_{F})]$ defines an isomorphism of $H^*(B)$-algebras from ${\mathfrak R}/{\mathfrak I} \rightarrow H^*(E(X))$. Here $V_F$ denotes the connected $T$-stable submanifold $\Upsilon^{-1}(F)$ of $X$ corresponding to a face $F$ of $Q$ and $[E(V_F)]$ denotes the ${Poincar\acute{e}}$ dual of $E(V_{F}):=E\times_{T}V_{F}$ in $H^*(E(X))$. \end{thm} {\bf Proof:} By \cite[Corollary 7.8]{masudapan} it follows that $H^*(X)$ is a free $\mathbb{Z}$-module of rank $\chi(X)$ and that there exists $p_1,\ldots, p_m$ polynomials in $\mathbb{Z}[x_{F}: F ~\mbox{a face of} ~Q]$ such that $p_i([V_{F}])$ for $1\leq i\leq m$ form a $\mathbb{Z}$-module basis of $H^*(X)$. Since $E(V_{F})\mid_{X}=V_{F}$ for each face $F$ of $Q$, by the Leray-Hirsch theorem $P_i:=p_i([E(V_F)])$ for $1\leq i\leq m$ form a basis of $H^*(E(X))$ as an $H^*(B)$-module.

Recall from (\ref{associatediso}) that we have the isomorphism of line bundles $\displaystyle{\prod_{i=1}^d (E(L_i))^{\langle u, v_i\rangle} \simeq \pi^*(\xi_u)}$ over $E(X)$ (see Lemma \ref{1}, Remark \ref{hompol} for the definition of the line bundles $L_i$ on $X$). Since $c_1(E(L_i))=[E(V_i)]$ we see that the relation $(ii)$ holds in $H^*(E(X))$ by (\ref{chernassociatediso}).

Consider the classifying map $f:B\longrightarrow BT$ of the principal $T$-bundle $E\longrightarrow B$. Thus we have the map $\widetilde{f}:E(X)\longrightarrow ET\times_{T} X$ over $f$ since $E(X)$ is the pull back of $ET\times_{T}X$ under $f$. This induces the canonical maps of cohomology rings $\widetilde{f}^*: H_{T}^*(X)\longrightarrow H^*(E(X))$ over $f^*:H^*(BT)\longrightarrow H^*(B)$ giving a commuting square
\[\begin{array}{lllllll}
H_{T}^*(X)&\stackrel{\widetilde{f}^*}{\longrightarrow} & H^*(E(X))\\
~~~~\uparrow~{\pi}'^* &                            &~~~~\uparrow~\pi^*\\
H^*(BT)&\stackrel{f^*}{\longrightarrow} & H^*(B)
\end{array}\](see Remark \ref{12}).

Furthermore, the submanifold $ET\times_{T} V_{F}$ of $ET\times_T X$ pulls back to the submanifold $E(V_{F})$ of $E(X)$ under $f$. Thus the class $\tau_{F}:=[ET\times_{T} V_F]\in H_{T}^*(X)$ maps to the class $\tau'_{F}:=[E(V_{F})]$ in the cohomology ring $H^*(E(X))$. This in particular implies that the element $\displaystyle{\tau_{G}\tau_{H}-\tau_{G\vee H}\sum_{E\in G\cap H}\tau_{E}}$ maps to $\displaystyle{\tau'_{G}\tau'_{H}-\tau'_{G\vee H}\sum_{E\in G\cap H}\tau'_{E}}$ in $H^*(E(X))$. However, by \cite[Theorem 7.7]{masudapan}, $\displaystyle{\tau_{G}\tau_{H}-\tau_{G\vee H}\sum_{E\in G\cap H}\tau_{E}}=0$ in $H_{T}^*(X)$. Hence the relation $(i)$ holds in $H^*(E(X))$.  Thus $\Phi_1$ induces a well defined map from $\mathfrak{R}/\mathfrak{I}\longrightarrow H^*(E(X))$. On the other hand \cite[Theorem 7.7, Corollary 7.8]{masudapan} imply that as an $H^*(BT)$-algebra, the ring $H^*_{T}(X)$ has the presentation $\mathfrak{R}'/\mathfrak{I}'$, where $\mathfrak{R}':=H^*(BT)[x_{F}: F~ \mbox{a face of}~ Q]$ and $\mathfrak{I}'$ is the ideal in $\mathfrak{R}'$ generated by the relations $(i)$ above and the relations $(ii)'$ $\displaystyle{\sum_{i=1}^d \langle u,v_i\rangle x_{Q_i}}-u$ for $u\in Hom(T,S^1)=H^2(BT)$. This further implies that $\mathfrak{R}/\mathfrak{I}$ is isomorphic to the ring $\displaystyle{H_{T}^*(X)\otimes _{H^*(BT)} H^*(B)}$ which is a free $H^*(B)$-module of rank $\chi(X)$ as in Lemma \ref{4} above. Here $H^*(B)$ is a $H^*(BT)$-module by the map $f^*$ which sends $u\in Hom(T,S^1)$ to the class $c_1(\xi_u)\in H^*(B)$.

Since $H^*(B)$ is a finitely generated abelian group the proof follows by the arguments similar to the proof of Theorem \ref{2}.  $\Box$

\begin{ex}{\em (see \cite[Example 3.2, Example 5.8]{masudapan} and \cite[Example $7.4.36$]{bp}) Let $X=S^4$ be the $4$-sphere identified with the following subset \[\{(z_1,z_2,y)\in \mathbb{C}^2\times \mathbb{R}: |z_1|^2+|z_2|^2+|y|^2=1\}.\] Define a $T=S^1\times S^1$-action on $X$ given by $(t_1,t_2)\cdot (z_1,z_2,y)=(t_1z_1,t_2z_2,y)$. The $T$ action on $X$ is locally standard with $X/T$ homeomorphic to \[Q=\{(x_1,x_2,y)\in \mathbb{R}^{3} :x_1^2+x_2^2+y^2=1, x_1\geq 0, x_2\geq 0\}.\] It has $2$ characteristic submanifolds $z_1=0$ and $z_2=0$. The intersection of the two characteristic submanifolds is disconnected and it is the union of the two $T$-fixed points $(0,0,1)$ and  $(0,0,-1)$. The circle subgroup of $T$ which fixes $\{z_1=0\}$ is given by $\{(t,1)~:~t \in S^1\}$ which corresponds to $e_1 \in \text{ Hom }(S^1, T)\cong \Z^2$. Similarly the circle subgroup fixing $\{z_2=0\}$ is given by $\{(1,t)~:~t \in S^1\}$ which corresponds to $e_2 \in \text{ Hom }(S^1, T)\cong \Z^2$. Here $e_1, e_2$ are the standard basis of $\Z^2$. Here the orbit space $Q$ is a $2$-ball with two $0$-faces denoted by $a$ and $b$ respectively and two $1$-faces denoted by $G$ and $H$ respectively. Thus the orbit space is not a homology polytope, but is a face-acyclic manifold with corners.

  Let $B=\mathbb{C}\mathbb{P}^1$ and $E\longrightarrow B$ denote the principal $T$-bundle associated to the direct sum of the line bundles $\co \oplus \co(1)$ where $\co$ denotes the trivial line bundle and $\co(1)$ denotes the tautological line bundle on $\mathbb{C}\mathbb{P}^1$. Consider the associated $S^4$ bundle $E(S^4)$ over $B$. By Theorem  \ref{cohombuntorus}, $H^*(E(S^4))$ has the presentation $\mathfrak{R}/\mathfrak{I}$ where $\mathfrak{R}=H^*(\mathbb{C}\mathbb{P}^1)[x_G, x_H,x_a,x_b]$ with $x_a$ and $x_b$ are of degree $2$ and $x_{G}$ and $x_{H}$ are of degree $4$ and $\mathcal{I}$ is the ideal in $\mathfrak{R}$ generated by the following two relations 
  $(i)$ $x_G\cdot x_{H}-x_a-x_b;~~ x_ax_b$ ~~$(ii)$ $x_{G}-c_1(\co)=x_{G}-0=x_{G}; ~~x_{H}-c_1(\co(1))$. }
\end{ex}

\subsection{$K$-ring of a torus manifold bundle}
Since $H^{odd}(X)=0$ the Atiyah Hirzebruch spectral sequence with $E_2^{p,q}=H^p(X; K^q(pt))$ collapses at the $E_2$ term and converges to $K^{p+q}(X)$ (see \cite[p. 208 ]{atiyah_hirzebruch_adams_shepherd_1972}). Moreover, since $H^*(X)$ is free abelian of rank $\chi(X)$ by \cite[p. 209]{atiyah_hirzebruch_adams_shepherd_1972} we have $K^r(X)=0$ when $r$ is odd and $K^r(X)\simeq \mathbb{Z}^m$ when $r$ is even. Here $m=\chi(X)$ is also equal to the number of vertices of $Q$.  In particular, $K^0(X)$ is free abelian of rank $m$.

Let $E\longrightarrow B$ be a principal $T$-bundle and $E(X):=E\times_{T} X$ the associated bundle over a base $B$ having the homotopy type of a finite CW complex.
  
Let $\mathfrak{S}:=K^*(B)[x_{F}: F~ \mbox{a face of}~ Q]$ and $\mathfrak{J}$ denote the ideal in $\mathfrak{S}$ defined by the following relations:

$(i)$ $\displaystyle{x_Gx_{H}-x_{G\vee H}\sum_{E\in G\cap H} x_{E}};$  
  
$(ii)$ $\displaystyle{\prod_{i: \langle u,v_i\rangle>0} (1-x_{Q_i}})^{\langle u,v_i\rangle} -[\xi_{u}]\prod_{i: \langle u,v_i\rangle<0} (1-x_{Q_i}) ^{-\langle u,v_i\rangle}$ for $u\in Hom(T,S^1)$.

We have the following conjecture on $K^*(E(X))$ as a $K^*(B)$-algebra. When $B=pt$ this shall give a presentation of the $K$-ring of $X$ which will generalize Sankaran's result stated in Theorem \ref{kring}. For arbitrary $B$ this shall generalize our Theorem \ref{6} proved above.

\begin{conj}\label{kringtorusbun}
The ring $K^*(E(X))$ is a free $K^*(B)$ module of rank $m=\chi(X)$ and is isomorphic to $\mathfrak{S}/\mathfrak{J}$.
\end{conj}

\begin{rmk}{\em 
  The difficulty in this case is because the cohomology is not generated in degree $2$ (see \cite[Example 4.10]{masudapan}), we cannot find canonical complex line bundles whose classes generate the $K$-ring as in \cite[Section 3]{param}.

  On the other hand it may be useful to define the analogue of the $K$-theoretic face ring $\mathcal{K}'(Q)$ when $Q$ is a nice manifold with corners so that when $Q$ is a homology polytope it agrees with $\mathcal{K}(Q)$ (see  Definition \ref{14}). One can then check whether $\mathcal{K}'(Q)$ has the structure of a free $RT$-module of rank $m$ generalizing the  Proposition \ref{8} above.

Consider the fibration \(E(X) \longrightarrow B\) where \(X\) is as above. When $B$ is a path connected, finite-dimensional
  CW-complex then by \cite[Theorem \(9.22\)]{daviskirk}, there exists
  a cohomology spectral sequence with
  \(E^{p,q}_2=H^p(B, K^q(X)) \Rightarrow K^{p+q}(E(X))\). Since
  $K^q(X)=0$ for $q$ odd this spectral sequence collapses at the $E_2$ term. We wonder if this gives enough information to deduce the structure of
  \(K^*(E(X))\) as a \(K^*(B)\)-module.}
\end{rmk}


\begin{thebibliography}{9}

\bibitem{atiyah_hirzebruch_adams_shepherd_1972}
  \uppercase{Atiyah, M. F. --- Hirzebruch, F. --- Adams,  J. F. --- Shepherd, G. C.:}
		\newblock {\em Vector bundles and homogeneous spaces}, pages 196--222.
		\newblock London Mathematical Society Lecture Note Series. Cambridge University
		Press, 1972.
		
		\bibitem{AM}
		\uppercase{ Atiyah, M. F. --- Macdonald, I. G.:}
		\newblock {\em Introduction to commutative algebra}.
		\newblock Addison-Wesley Publishing Co., Reading, Mass.-London-Don Mills, Ont.
		1969.
		
		\bibitem{baggio2007equivariant}
		\uppercase{Baggio, S.:}
		\newblock \textit{Equivariant {$K$}-theory of smooth toric varieties,}
		\newblock Tohoku Mathematical Journal Second Series \textbf{59}(2) (2007), 203--231.
		
		\bibitem{bredon}
		\uppercase{Bredon, G. E.:}
		\newblock {\em Introduction to compact transformation groups}.
		\newblock Academic Press, New York-London, 1972.
		\newblock Pure and Applied Mathematics, Vol. \textbf{46}.
		
		\bibitem{bp}
		\uppercase{Buchstaber, V. M. --- Panov, T. E.:}
		\newblock {\em Toric topology}, volume 204 of {\em Mathematical Surveys and
			Monographs}.
		\newblock American Mathematical Society, Providence, RI, 2015.
		
		\bibitem{dan}
		\uppercase{Danilov, V. I.:}
		\newblock \textit{The geometry of toric varieties}
		\newblock Uspekhi Mat. Nauk, \textbf{33}(2(200)) :85--134, 247, 1978.
		
		\bibitem{dku}
		\uppercase{Dasgupta, J.--- Khan, B. --- Uma, V.:}
		\newblock \textit{Equivariant {$K$}-ring of quasitoric manifolds,}
		\newblock  ArXiv e-prints; arXiv:1805.11373 [math.AT], 2018.
		
		\bibitem{daviskirk}
		\uppercase{Davis, J. F. --- Kirk, P.:}
		\newblock {\em Lecture notes in algebraic topology}, volume~\textbf{35} of {\em Graduate
			Studies in Mathematics}.
		\newblock American Mathematical Society, Providence, RI, 2001.
		
		\bibitem{davisjanu}
		\uppercase{Davis, M. W. --- Januszkiewicz, T.:}
		\newblock \textit{Convex polytopes, {C}oxeter orbifolds and torus actions,}
		\newblock  Duke Math. J. \textbf{62}(2) (1991), 417--451.
		
		\bibitem{ful}
		\uppercase{Fulton, W.:}
		\newblock {\em Introduction to toric varieties}, volume 131 of {\em Annals of
			Mathematics Studies}.
		\newblock Princeton University Press, Princeton, NJ, 1993.
		\newblock The William H. Roever Lectures in Geometry.
		
		\bibitem{moment}
		\uppercase{ Guillemin, V. --- Ginzburg, V. --- Karshon, Y.:}
		\newblock {\em Moment maps, cobordisms, and {H}amiltonian group actions},
		volume~\textbf{98} of {\em Mathematical Surveys and Monographs}.
		\newblock American Mathematical Society, Providence, RI, 2002.
		\newblock Appendix J by Maxim Braverman.
		
		\bibitem{top}
		\uppercase{Hatcher, A.:}
		\newblock {\em Algebraic topology}.
		\newblock Cambridge University Press, Cambridge, 2002.
		
		\bibitem{hatcher2003vector}
		\uppercase{Hatcher, A.:}
		\newblock \textit{Vector bundles and {$K$}-theory,}
		\newblock In Internet under http://www. math. cornell. edu/\~{} hatcher,
		2003.
		
		\bibitem{hm}
		\uppercase{Hattori, A. --- Masuda, M.:}
		\newblock \textit{Theory of multi-fans,}
		\newblock  Osaka J. Math. \textbf{40}(1) (2003), 1--68.
		
		\bibitem{huneke}
		\uppercase{Huneke, C. --- Swanson, I.:}
		\newblock {\em Integral closure of ideals, rings, and modules}, volume \textbf{336} of
		{\em London Mathematical Society Lecture Note Series}.
		\newblock Cambridge University Press, Cambridge, 2006.
		
		\bibitem{karoubi}
		\uppercase{Karoubi, M.:}
		\newblock {\em {$K$}-theory}.
		\newblock Classics in Mathematics. Springer-Verlag, Berlin, 2008.
		\newblock An introduction, Reprint of the 1978 edition, With a new postface by
		the author and a list of errata.
		
		\bibitem{masudapan}
		\uppercase{ Masuda, M. --- Panov, T.}
		\newblock \textit{On the cohomology of torus manifolds,}
		\newblock  Osaka J. Math. \textbf{43} (2006), 711--746.
		
		\bibitem{gmukherjee}
		\uppercase{Mukherjee, G., editor:}
		\newblock {\em Transformation groups}
		\newblock Hindustan Book Agency, New Delhi, 2005.
		\newblock Symplectic torus actions and toric manifolds, With contributions by
		Chris Allday, Mikiya Masuda and P. Sankaran.
		
		\bibitem{param}
		\uppercase{Sankaran, P.:}
		\newblock \textit{{$K$}-rings of smooth complete toric varieties and related spaces,}
		\newblock  Tohoku Math. J. \textbf{60} (2008), 459--469.
		
		\bibitem{paramuma}
		\uppercase{Sankaran, P. --- Uma, V.:}
		\newblock \textit{Cohomology of toric bundles,}
		\newblock Comment. Math. Helv. \textbf{78}(4) (2003), 540--554.
		
		\bibitem{MR2313027}
		\uppercase{Sankaran, P. --- Uma, V.:}
		\newblock \textit{{$K$}-theory of quasitoric manifolds,}
		\newblock Osaka J. Math. \textbf{44}(1) (2007), 71--89.
		
		\bibitem{suyama}
		\uppercase{Suyama, Y.:}
		\newblock \textit{Examples of smooth compact toric varieties that are not quasitoric
			manifolds,}
		\newblock Algebr. Geom. Topol. \textbf{14}(5) (2014), 3097--3106.
		
		\bibitem{uma}
		\uppercase{Uma, V.:}
		\newblock \textit{{$K$}-theory of torus manifolds,}
		\newblock Toric topology, Contemp. Math. \textbf{460} (2008), 85--389.
		
		\bibitem{vezzosi2003higher}
		\uppercase{Vezzosi, G. --- Vistoli, A.:}
		\newblock \textit{Higher algebraic {$K$}-theory for actions of diagonalizable groups,}
		\newblock Inventiones mathematicae \textbf{153}(1) (2003), 1--44.


\end{thebibliography}
\end{document}